\documentclass[11pt]{article}
\newcommand{\bbR}{{\rm I\hspace{-0.7mm}R}}

\usepackage{latexsym}

\begin{document}

\centerline{\bf \large On A General Theorem For Additive L\'evy Processes }

\vskip .1in

\centerline{\tt By MING YANG}

\vskip .1in

\indent{\footnotesize

We prove a new theorem and show that this theorem implies
Theorem 1.1, Theorem 1.5 (with no restrictions), Theorem 2.1 of Khoshnevisan, Xiao and Zhong [1]
and Theorem 2.2 of Khoshnevisan and Xiao [2]. }

\vskip .1in
\noindent
{\small \noindent 2000 {\it Mathematics Subject Classification}. Primary
60G60, 60G51; secondary 60G17.
  
}

\vskip .15in
\centerline{\textsc{1. Introduction}}
\vskip .15in

Let
$X^1_{t_1},~X^2_{t_2},\cdots,X^N_{t_N}$ be $N$ independent L\'evy
processes in $\bbR^d$ with their respective L\'evy exponents 
$\Psi_j,~j=1,2,\cdots,N$. The random field
$$X_t=X^1_{t_1}+X^2_{t_2}+\cdots+X^N_{t_N},~~~~~t=(t_1,t_2,\cdots,t_N)\in\bbR^N_+$$
is called the additive
L\'evy process. Let $\lambda_d$ denote Lebesgue measure in $\bbR^d.$
Define $E_1+E_2=\{x+y:x\in E_1,~y\in E_2\}$ 
for any two sets $E_1,~E_2$ of $\bbR^d.$

\vskip 0.15in
\noindent
{\bf Theorem 1.1}
{\it
Let $(X;~\Psi_1,\cdots,\Psi_N)$ be an additive L\'evy process in $\bbR^d$ 
and $(Y;~\Psi'_1,\cdots,\Psi'_n)$ be a second
additive L\'evy process in $\bbR^d$ independent of $X$.
Then for any
$G\in\mathcal{B}(\bbR^n_+)\backslash\{\emptyset\}$,} 
$$E\{\lambda_d(X(\bbR^N_+)+Y(G))\}>0\Longleftrightarrow$$
$$\int_{\bbR^d}\left(
\int_{\bbR^n_+}\int_{\bbR^n_+}
e^{-\sum^n_{i=1}|t_i-s_i|\Psi'_i(\mbox{sgn}(t_i-s_i)\xi)}\mu(ds)\mu(dt)\right)
\prod_{j=1}^N\mbox{Re}\left(\frac{1}{1+\Psi_j(\xi)}\right)
d\xi<\infty\eqno(1.1)$$
{\it for some probability measure $\mu$ on $G$.} 

\vskip .15in

If we take $G=\{0\}$, 
we obtain
$$E\{\lambda_d(X(\bbR^N_+))\}>0\Longleftrightarrow
\int_{\bbR^d}
\prod_{j=1}^N\mbox{Re}\left(\frac{1}{1+\Psi_j(\xi)}\right)
d\xi<\infty.\eqno(1.2)$$
This is 
Theorem 1.1 of Khoshnevisan, Xiao and Zhong [1] with no assumptions.
If we let $X=0$ we obtain Theorem 2.1 of
Khoshnevisan, Xiao and Zhong [1]. If we take a standard additive
$\alpha$-stable process $S^{\alpha}$ 
from $\bbR^d_+$ to $\bbR^d$ with $\alpha\in(0,1)$ to be $X$,
we obtain Theorem 2.2 of Khoshnevisan and Xiao [2].
If we consider
a deterministic additive L\'evy process from $\bbR^d_+$ to $\bbR^d$ as $Y$, 
we obtain 

\vskip .15in
{\it
For any $F\in\mathcal{B}(\bbR^d)\backslash\{\emptyset\},$}
$$E\{\lambda_d(X(\bbR^N_+)+F)\}>0\Longleftrightarrow
\int_{\bbR^d}|\hat{\mu}(\xi)|^2
\prod_{j=1}^N\mbox{Re}\left(\frac{1}{1+\Psi_j(\xi)}\right)
d\xi<\infty\eqno(1.3)$$
{\it for some probability measure $\mu$ on $F$, where
$\hat{\mu}(\xi)=\int_{\bbR^d}e^{i\xi\cdot x}\mu(dx),~\xi\in\bbR^d.$} 

\vskip .15in
\noindent
(The proof is given in the next section.) This is 
Theorem 1.5 of Khoshnevisan, Xiao and Zhong [1] without extra conditions.

\vskip .15in
\centerline{\textsc{2. Proof of Theorem 1.1}}
\vskip .15in

Our proof is based on a modification of 
the proof of Theorem 2.2 of  
Khoshnevisan and Xiao [2]. 
Let us lay out some groundwork first.

Let $z_1,\cdots,z_N$ be $N$ complex numbers such that
Re$(z_j)\ge0.$ There are $2^N$ 
different permutations to write down a
partial conjugate of
the vector
$(z_1,\cdots,z_N).$
For example, 
$(z_1,\cdots,z_N)$ (itself),
$(\overline{z_1},\cdots,z_N),$ 
$(z_1,\cdots,\overline{z_N}),$ 
$(\overline{z_1},\overline{z_2},\cdots,z_N),$
$(\overline{z_1},\overline{z_2},\cdots,\overline{z_N}),$  
and so on. Let 
$(z_1^{\pm},z_2^{\pm},\cdots,z_N^{\pm})$ denote the generic partially conjugated vector.
Then, we have
$$\sum_{\pm}
\prod_{j=1}^N\frac{1}{1+z_j^{\pm}}=2^N\prod_{j=1}^N\mbox{Re}
\left(\frac{1}{1+z_j}\right)>0.\eqno(2.1)$$ 
Identity (2.1) can be proved by induction.
It follows immediately from (2.1) that
$$\sum_{\pm}
\prod_{j=1}^N\frac{1}{1+\Psi_j^{\pm}(\xi)}=2^N\prod_{j=1}^N\mbox{Re}
\left(\frac{1}{1+\Psi_j(\xi)}\right)>0.\eqno(2.2)$$
Let $X^{\pm}$ be the additive L\'evy process with L\'evy exponent
$(\Psi_1^{\pm},\cdots,\Psi^{\pm}_N)$. In other words, $X^{\pm}=\pm
X^1\pm X^2\cdots\pm X^N.$ Here, $\pm$ has the true meaning of $+$ or $-.$
Define 
$$Z^{\pm}=Y+X^{\pm}.$$
Here we require the $2^N$ additive L\'evy processes $Z^{\pm}$ to be totally independent.
In other words,
we have $2^N$ independent copies of $Y$ (with the same notation though) and
$2^N$ independent additive L\'evy processes $X^{\pm}$, independent of the $Y$'s as well.
For the sake of convenience, we index the
$Z^{\pm}$ as
$Z^1,Z^2,\cdots,Z^{2^N}$ with
$Z^1=Y+X.$ We define a super
additive L\'evy process
$$Z=(Z^1,Z^2,\cdots,Z^{2^N}).$$
Clearly, $Z$ is a $2^N(n+N)-$parameter additive L\'evy process taking values
in $\bbR^{2^Nd}.$

Let $P_{\lambda_{2^Nd}}$ and $E_{\lambda_{2^Nd}}$ be the sigma-finite measure and the corresponding
expectation with respect to $Z.$
The reader is referred to 
Sections 3, 4 of Khoshnevisan, Xiao and Zhong [1] for all developments about
$P_{\lambda_{2^Nd}}$ and $E_{\lambda_{2^Nd}}.$
Throughout, we only write
$P_{\lambda}$ and $E_{\lambda}$
rather than $P_{\lambda_{2^Nd}}$ and $E_{\lambda_{2^Nd}},$ respectively.
We also introduce the $2^N(n+N)-$parameter process
$\mathcal{M}^{A,f,\mu}_s$ based on $Z,$ the same thing as $\mathcal{M}^{A}_{\mu}f(s)$ 
in
Khoshnevisan, Xiao and Zhong [1]. One of our proof techniques is to manipulate the four 
parameters $A,f,\mu,s.$

\vskip 0.15in
\noindent
{\bf Proof of Theorem 1.1}~~Since the
direction $\Longleftarrow$ is a special case
of Theorem 2.1 of
Khoshnevisan, Xiao and Zhong [1], the direction 
$\Longrightarrow$ is something we have to prove. Suffice it to show that
$$E\{\lambda_d(X([0,l]^N)+Y(G))\}>0\Longrightarrow$$
$$\int_{\bbR^d}\left(
\int_{\bbR^n_+}\int_{\bbR^n_+}
e^{-\sum^n_{i=1}|t_i-s_i|\Psi'_i(\mbox{sgn}(t_i-s_i)\xi)}\mu(ds)\mu(dt)\right)
\prod_{j=1}^N\mbox{Re}\left(\frac{1}{1+\Psi_j(\xi)}\right)
d\xi<\infty$$
for some probability measure $\mu$ on $G$, where
$G$ is compact and $l\in(0,\infty).$

By Proposition 10.3 and Theorem 2.1 of Khoshnevisan, Xiao and Zhong [1], it is always true that
$E\{\lambda_d(X'([0,l]^N)+Y'(G))\}<\infty$ for any 
processes $X',~Y'$ and compact $G$. We  
separate the two cases whether $G$ has 
positive Lebesgue measure.

\vskip .15in
\noindent
{\bf Case 1} $\lambda_n(G)>0.$~~~
In this case, by Proposition 10.3 and Theorem 2.1 of Khoshnevisan, Xiao and Zhong [1], 
$$E\{\lambda_d(Z^i(G\times[0,l]^N))\}\approx
E\{\lambda_d(Z^j(G\times[0,l]^N))\}\eqno(2.3)$$
for any two $Z^i,~Z^j$, $\approx$ depending on $\lambda_n(G)$ and $l$.
[Here, the symbol $\approx$ means that there is a constant $c\in(0,\infty)$
such that $c^{-1}f_1\le f_2\le cf_1$ for two nonnegative functions or quantities $f_1$ and $f_2.$]

Let $G^{\delta}$ be the
closed $\delta$-enlargement of $G$ for $\delta>0,$ that is,
the smallest compact set such that
for each point $s=(s_1,\cdots,s_n)\in G,$
$[s_1,s_1+\delta]\times\cdots\times [s_n,s_n+\delta]\subset G^{\delta}.$
Let $|x|=(x\cdot x)^{1/2}$
and let $B(x,r)$ be the closed ball of radius $r$ with center at $x$.
Define $B^{\delta}=(B(0,\delta))^{2^N}=
(B(0,\delta),\cdots,B(0,\delta)).$

By the definition of 
$P_{\lambda},$ also thanks to
Fubini's theorem, independence and the fact that

\noindent
$-B(0,r)=B(0,r),$ 
we have
\begin{eqnarray*}
&&P_{\lambda}
\{Z((G^{\delta}\times[0,l]^N)^{2^N})\bigcap
B^{\delta}\neq\emptyset\}\\
&&=
\int_{\bbR^{2^Nd}}P\left\{(x+Z((G^{\delta}\times[0,l]^N)^{2^N}))\bigcap
B^{\delta}\neq\emptyset\right\}dx\\
&&=
\int_{\bbR^{2^Nd}}P\left\{
Z^1(G^{\delta}\times[0,l]^N)\bigcap
B(x_1,\delta)\neq\emptyset,\cdots,
Z^{2^N}(G^{\delta}\times[0,l]^N)\bigcap
B(x_{2^N},\delta)\neq\emptyset\right\}\\
&&
~~~\cdot dx_1\cdots dx_{2^N}\\
&&=
\prod_{i=1}^{2^N}\int_{\bbR^{d}}P\left\{
Z^i(G^{\delta}\times[0,l]^N)\bigcap
B(x,\delta)\neq\emptyset\right\}dx\\
&&=
\prod_{i=1}^{2^N}
E\left\{
\lambda_d(Z^i(G^{\delta}\times[0,l]^N)+
B(0,\delta))\right\}\\
&&\rightarrow
\prod_{i=1}^{2^N}
E\left\{
\lambda_d\left(\overline{Z^i(G\times[0,l]^N)}\right)\right\}
\end{eqnarray*}
downwards as $\delta\rightarrow0$.
By (2.3),
$$
\prod_{i=1}^{2^N}
E\left\{
\lambda_d\left(\overline{Z^i(G\times[0,l]^N)}\right)\right\}
\ge
c_1\left[E\{\lambda_d(X([0,l]^N)+Y(G))\}\right]^{2^N}>0$$
for some constant $c_1\in(0,\infty),$ and for all $\delta>0,$
$$P_{\lambda}
\{Z((G^{\delta}\times[0,l]^N)^{2^N})\bigcap
B^{\delta}\neq\emptyset\}\in(0,\infty).$$

We add a cemetery point
$\Delta\notin\bbR^{n}_+$ to $\bbR^n_+$ to construct
a measurable 
map $T^{\delta}$ (random variable) from $\Omega$ to ${\bf Q}^n_+\cup\{\Delta\},$
where ${\bf Q}$ stands for rational as always.
$T^{\delta}$ is defined as follows.
$T^{\delta}\neq\Delta$ if and only if 
$T^{\delta}\in{\bf Q}^n_+\bigcap(0,\infty)^n\bigcap G^{\delta}$
and
there exist
$t_i\in (0,l]^N\bigcap
{\bf Q}^N_+,$
$1\le i\le 2^N,$ such that
$|Z^{i}_{(T^{\delta},t_i)}|\le\delta.$
This can always be done.
We have
\begin{eqnarray*}
&&P_{\lambda}\left\{T^{\delta}\neq\Delta\right\}\\
&&=P_{\lambda}
\left\{Z((G^{\delta}\times[0,l]^N)^{2^N})\bigcap
B^{\delta}\neq\emptyset\right\}\\
&&=P_{\lambda}
\left\{Z^i(G^{\delta}\times[0,l]^N)\bigcap
B(0,\delta)\neq\emptyset,~1\le i\le2^N\right\}\in(0,\infty).
\end{eqnarray*}
There is therefore
a probability measure $\mu^{\delta}$ in $\bbR^n_+$ supported
on $G^{\delta}$ given by
$$\mu^{\delta}(\bullet)=
\frac{P_{\lambda}\{T^{\delta}\in\bullet,~T^{\delta}\neq\Delta\}}
{P_{\lambda}\{T^{\delta}\neq\Delta\}}.\eqno(2.4)$$

By Lemma 4.2 of Khoshnevisan, Xiao and Zhong [1], for any
${\bf A},~{\bf f},~\nu$ [$\nu$ probability measure in $\bbR_+^{2^N(n+N)}$] and
${\bf s}\in\bbR_+^{2^N(n+N)}$,
we have
$$\mathcal{M}^{{\bf A},{\bf f},\nu}_{{\bf s}}
\ge
\int_{{\bf t}\succeq_{{\bf A}} {\bf s}}
P_{{\bf t}-{\bf s}}{\bf f}(Z_{{\bf s}})\nu(d{\bf t})
,~~~P_{\lambda}-a.s.\eqno(2.5)$$
See (3.1) and Lemma 3.1 of Khoshnevisan, Xiao and Zhong [1] for
the definition of the operator
$P_{{\bf t}}{\bf f}.$
Thus,
$$\mathcal{M}^{{\bf A},{\bf f},\nu}_{{\bf s}}
\ge
\int_{{\bf t}\succeq_{{\bf A}} {\bf s}}
P_{{\bf t}-{\bf s}}{\bf f}(Z_{{\bf s}})\nu(d{\bf t})
\cdot
1_{\{|Z^{i}_{{\bf s}^i}|\le\delta ,~1\le i\le2^N \}}
,~~~P_{\lambda}-a.s.\eqno(2.6)$$
where ${\bf s}=({\bf s}^1,\cdots,{\bf s}^{2^N}),~{\bf s}^i\in\bbR_+^{n+N}.$ Our next step is to
make $2^N$ different combinations of
${\bf A},{\bf f},\nu, {\bf s}$ so that we have 
$2^N$ different inequalities of (2.6), and then
we sum them up to see what happens to the right-hand side.
Here, the $P_{\lambda}$-null set in (2.5) depends on ${\bf s}$. Thus, if
${\bf s}$ is random and if we wish that (2.5) holds 
uniformly in $\omega$, one way is to require
${\bf s}$ to take rational points only.

For $\varepsilon>0,$
define
$f_{\varepsilon}(x)=(2\pi\varepsilon)^{-d/2}e^{-|x|^2/2\varepsilon},~x\in\bbR^d.$
Let $u,~v\in\bbR^n_+$ and $s,~t\in\bbR^N_+.$ In this paper,
we have only to consider the partial orders on
$\bbR^n_+.$ So, if $\pi$ is a partial order on $\bbR^n_+,$
the corresponding partial order $A$ on $\bbR^{n+N}_+$ is defined by
$(u,s)\preceq_A(v,t)\Longleftrightarrow
u\preceq_{\pi} v,~s_j\le t_j,~1\le j\le N.$
Let $\kappa(dt)=e^{-\sum^N_{j=1}t_j}dt,~t=(t_1,\cdots,t_N)\in\bbR^N_+$ 
and let $\mu$ be
a probability measure
in $\bbR^n_+.$ Take
$${\bf s}=(0,\cdots0,
(v,s)_i,0\cdots,0),~~s\in\bbR^N_+,~~v\in\bbR^n_+,\eqno(2.7)$$
$$\nu=\mu\otimes\kappa\otimes\mu_0,\eqno(2.8)$$
$${\bf f}({\bf x})=f_{\varepsilon}(x)e^{-\frac{|x'|^2}{2}}.\eqno(2.9)$$
The notation in (2.7) is clear. In (2.8), $\mu_0$ is the point mass at $0$ acting
on all irrelevant time parameters other than $(v,s)_i.$
Note that
for any L\'evy exponent $\Psi'',$
$$\int_{\bbR^m_+}
e^{-\sum^m_{i=1}u_i\Psi''_i(\xi)}\mu_0(du)=1.$$
In (2.9), $x$ corresponds to $Z^i$ and
$x'$ is for all other irrelevant components in the value space of $Z$.
For ${\bf A}$, we take
${\bf A}=(A_0,A)$. Here, $A_0$ is always the componentwise natural order, i.e.,
$(a_1,\cdots,a_m)\preceq_{A_0}
(b_1,\cdots,b_m)\Longleftrightarrow
a_k\le b_k,~1\le k\le m.$
Of course, $A$ is for $(v,s)_i$ while $A_0$ takes place elsewhere.  
By Lemma 3.1 of
Khoshnevisan, Xiao and Zhong [1], we have
\begin{eqnarray*}
&&
\int_{{\bf t}\succeq_{{\bf A}} {\bf s}}
P_{{\bf t}-{\bf s}}{\bf f}(0)\nu(d{\bf t})\\
&&
=c_2e^{-\sum^N_{j=1}s_j}
\int_{u\succeq_{\pi} v}\int_{\bbR^d}
e^{-\frac{\varepsilon^2}{2}|\xi|^2}
e^{-\sum^n_{i=1}|u_i-v_i|\Psi'_i(\mbox{sgn}(u_i-v_i)\xi)}
\prod_{j=1}^N\frac{1}{1+\Psi_j^{\pm}(\xi)}d\xi\mu(du),
\end{eqnarray*}
where $c_2\in(0,\infty)$ is a constant which involves
$\pi,~d,~N,~n,$ as well as a $(2^N-1)d$
power of
$\int_{-\infty}^{\infty}e^{-\frac{y^2}{2}}dy$, but is identical for our 
all choices of $({\bf A},{\bf f},\nu,{\bf s}).$
Here, it is valid to interchange the order of integration owing to the
term $e^{-\frac{\varepsilon^2}{2}|\xi|^2}\cdot
e^{-\frac{1}{2}|\xi'|^2}.$
Therefore, by (2.2),
\begin{eqnarray*}
&&
\sum_{i=1}^{2^N}
\int_{{\bf t}\succeq_{{\bf A}} {\bf s}}
P_{{\bf t}-{\bf s}}{\bf f}(0)\nu(d{\bf t})\\
&&
=c_3e^{-\sum^N_{j=1}s_j}
\int_{u\succeq_{\pi} v}\int_{\bbR^d}
e^{-\frac{\varepsilon^2}{2}|\xi|^2}
e^{-\sum^n_{i=1}|u_i-v_i|\Psi'_i(\mbox{sgn}(u_i-v_i)\xi)}
\prod_{j=1}^N\mbox{Re}\left(\frac{1}{1+\Psi_j(\xi)}\right)d\xi\mu(du).
\end{eqnarray*}
($c_3=2^Nc_2.$)
On the other hand,
By Lemma 4.2 of
Khoshnevisan, Xiao and Zhong [1], 
$$E_{\lambda}
\left(
\sup_{\theta\in{\bf Q}^{2^N(n+N)}_+}\mathcal{M}^{{\bf A},{\bf f},\nu}_{\theta}
\right)^2
\le c_4
\int_{\bbR^d}e^{-\frac{\varepsilon^2}{2}|\xi|^2}Q_{\mu}(\xi)d\xi,\eqno(2.10)$$
where
$$Q_{\mu}(\xi)
=\int_{\bbR^n_+}\int_{\bbR^n_+}
e^{-\sum^n_{i=1}|t_i-s_i|\Psi'_i(\mbox{sgn}(t_i-s_i)\xi)}\mu(ds)\mu(dt)\cdot
\prod_{j=1}^N\mbox{Re}\left(\frac{1}{1+\Psi_j(\xi)}\right).$$
Here, $c_4\in(0,\infty)$ is a constant similar to $c_2.$ To justify (2.10),
we notice that
$e^{-\varepsilon^2|\xi|^2}<
e^{-\frac{\varepsilon^2}{2}|\xi|^2}$, 
$Q_{\mu}(\xi)\in(0,1],$ and for all $X^{\pm},$
$$\prod_{j=1}^N\mbox{Re}
\left(\frac{1}{1+\Psi_j^{\pm}(\xi)}\right)=
\prod_{j=1}^N\mbox{Re}
\left(\frac{1}{1+\Psi_j(\xi)}\right),$$
and one final detail, the effect of the point mass
$\mu_0.$

The Lipschitz continuity of
${\bf f}$ is evident. Let
$D(\varepsilon)$ be the 
Lipschitz constant of ${\bf f}.$
By the definition of $P_{{\bf t}-{\bf s}}{\bf f}$,
$$D(\varepsilon)\delta
+\inf_{|z|\le\delta}P_{{\bf t}-{\bf s}}
{\bf f}(z)\ge P_{{\bf t}-{\bf s}}{\bf f}(0).$$
Since $\nu$ is a probability measure,
$$D(\varepsilon)\delta
+\int_{{\bf t}\succeq_{{\bf A}} {\bf s}}
\inf_{|z|\le\delta}P_{{\bf t}-{\bf s}}
{\bf f}(z)\nu(d{\bf t})
\ge \int_{{\bf t}\succeq_{{\bf A}} {\bf s}}
P_{{\bf t}-{\bf s}}{\bf f}(0)\nu(d{\bf t}).$$
If 
$|Z^{i}_{{\bf s}^i}|\le\delta ,~1\le i\le2^N,$
then
$|Z_{{\bf s}}|=|Z_{0,{\bf s}^i,0}|\le\delta$
and
$P_{{\bf t}-{\bf s}}{\bf f}(Z_{{\bf s}})\ge
\inf_{|z|\le\delta}P_{{\bf t}-{\bf s}}
{\bf f}(z).$ 
Note that
$\inf_{|z|\le\delta}P_{{\bf t}-{\bf s}}
{\bf f}(z)$ is a function of ${\bf t}$ independent of $\omega$ for each
fixed ${\bf s}.$
Thus,
\begin{eqnarray*}
&&\int_{{\bf t}\succeq_{{\bf A}} {\bf s}}
P_{{\bf t}-{\bf s}}{\bf f}(Z_{{\bf s}})\nu(d{\bf t})
\cdot
1_{\{|Z^{i}_{{\bf s}^i}|\le\delta ,~1\le i\le2^N \}}\\
&&\ge
\int_{{\bf t}\succeq_{{\bf A}} {\bf s}}
\inf_{|z|\le\delta}P_{{\bf t}-{\bf s}}
{\bf f}(z)\nu(d{\bf t})
\cdot
1_{\{|Z^{i}_{{\bf s}^i}|\le\delta ,~1\le i\le2^N \}}\\
&&\ge\left[
\int_{{\bf t}\succeq_{{\bf A}} {\bf s}}
P_{{\bf t}-{\bf s}}{\bf f}(0)\nu(d{\bf t})-
D(\varepsilon)\delta\right]
\cdot
1_{\{|Z^{i}_{{\bf s}^i}|\le\delta ,~1\le i\le2^N \}}
\end{eqnarray*}
and subsequently
\begin{eqnarray*}
&&\sum_{i=1}^{2^N}
\int_{{\bf t}\succeq_{{\bf A}} {\bf s}}
P_{{\bf t}-{\bf s}}{\bf f}(Z_{{\bf s}})\nu(d{\bf t})
\cdot
1_{\{|Z^{i}_{{\bf s}^i}|\le\delta ,~1\le i\le2^N \}}\\
&&\ge
\sum_{i=1}^{2^N}
\left[
\int_{{\bf t}\succeq_{{\bf A}} {\bf s}}
P_{{\bf t}-{\bf s}}{\bf f}(0)\nu(d{\bf t})-
D(\varepsilon)\delta\right]
\cdot
1_{\{|Z^{i}_{{\bf s}^i}|\le\delta ,~1\le i\le2^N \}}\\
&&=
[c_3e^{-\sum^N_{j=1}s_j}
\int_{u\succeq_{\pi} v}\int_{\bbR^d}
e^{-\frac{\varepsilon^2}{2}|\xi|^2}
e^{-\sum^n_{i=1}|u_i-v_i|\Psi'_i(\mbox{sgn}(u_i-v_i)\xi)}
\prod_{j=1}^N\mbox{Re}\left(\frac{1}{1+\Psi_j(\xi)}\right)d\xi\mu(du)\\
&&
~~~-2^N
D(\varepsilon)\delta]\cdot
1_{\{|Z^{i}_{{\bf s}^i}|\le\delta ,~1\le i\le2^N \}}.
\end{eqnarray*}
It follows from (2.6) and the definition of $T^{\delta}$ that
\begin{eqnarray*}
&&\sum_{i=1}^{2^N}
\sup_{\theta\in{\bf Q}^{2^N(n+N)}_+}\mathcal{M}^{{\bf A},{\bf f},\nu}_{\theta}\\
&&\ge
[
c_3e^{-lN}
\int_{u\succeq_{\pi} T^{\delta}}\int_{\bbR^d}
e^{-\frac{\varepsilon^2}{2}|\xi|^2}
e^{-\sum^n_{i=1}|u_i-T^{\delta}_i|\Psi'_i(\mbox{sgn}(u_i-T^{\delta}_i)\xi)}
\prod_{j=1}^N\mbox{Re}\left(\frac{1}{1+\Psi_j(\xi)}\right)d\xi\mu(du)\\
&&
~~~-2^N
D(\varepsilon)\delta]\cdot
1_{\{T^{\delta}\neq\Delta \}},~~~~~~~~~~~P_{\lambda}-a.s.
\end{eqnarray*}
We rewrite the preceding as 
$$
2^N
D(\varepsilon)\delta\cdot
1_{\{T^{\delta}\neq\Delta \}}+
\sum_{i=1}^{2^N}
\sup_{\theta\in{\bf Q}^{2^N(n+N)}_+}\mathcal{M}^{{\bf A},{\bf f},\nu}_{\theta}
\ge
c_3e^{-lN}
\int_{u\succeq_{\pi} T^{\delta}}\int_{\bbR^d}
e^{-\frac{\varepsilon^2}{2}|\xi|^2}$$
$$
e^{-\sum^n_{i=1}|u_i-T^{\delta}_i|\Psi'_i(\mbox{sgn}(u_i-T^{\delta}_i)\xi)}
\prod_{j=1}^N\mbox{Re}\left(\frac{1}{1+\Psi_j(\xi)}\right)d\xi\mu(du)
\cdot
1_{\{T^{\delta}\neq\Delta \}},~~~P_{\lambda}-a.s.\eqno(2.11)$$
Now, in (2.11) replace $\mu$ by $\mu^{\delta}$ and appeal to
the Cauchy-Schwarz inequality,
$$\left(
\sum^p_{i=1}x_i\right)^2\le
p\sum^p_{i=1}x_i^2\eqno(2.12)$$
for any $p$ real numbers $x_i,~i=1,\cdots p,$ to obtain
$$
2^{N+1}\delta^2
D^2(\varepsilon)\cdot
1_{\{T^{\delta}\neq\Delta \}}+
2^{N+1}\sum_{i=1}^{2^N}
\left(\sup_{\theta\in{\bf Q}^{2^N(n+N)}_+}\mathcal{M}^{{\bf A},{\bf f},\nu}_{\theta}\right)^2
\ge
c_3^2e^{-2lN}
\{\int_{u\succeq_{\pi} T^{\delta}}\int_{\bbR^d}
e^{-\frac{\varepsilon^2}{2}|\xi|^2}$$
$$
e^{-\sum^n_{i=1}|u_i-T^{\delta}_i|\Psi'_i(\mbox{sgn}(u_i-T^{\delta}_i)\xi)}
\prod_{j=1}^N\mbox{Re}\left(\frac{1}{1+\Psi_j(\xi)}\right)d\xi\mu^{\delta}(du)\}^2
\cdot
1_{\{T^{\delta}\neq\Delta \}},~~~P_{\lambda}-a.s.\eqno(2.13)$$
Taking $E_{\lambda}$-expectation on both sides of (2.13) followed by the Cauchy-Schwarz
inequality yields
$$
2^{N+1}\delta^2
D^2(\varepsilon)
P_{\lambda}\{T^{\delta}\neq\Delta \}+
2^{N+1}\sum_{i=1}^{2^N}
E_{\lambda}\left(\sup_{\theta\in{\bf Q}^{2^N(n+N)}_+}\mathcal{M}^{{\bf A},{\bf f},\nu}_{\theta}\right)^2
\ge
c_3^2e^{-2lN}
\{\int_{\bbR^n}\int_{u\succeq_{\pi} v}\int_{\bbR^d}$$
$$e^{-\frac{\varepsilon^2}{2}|\xi|^2}
e^{-\sum^n_{i=1}|u_i-v_i|\Psi'_i(\mbox{sgn}(u_i-v_i)\xi)}
\prod_{j=1}^N\mbox{Re}\left(\frac{1}{1+\Psi_j(\xi)}\right)d\xi\mu^{\delta}(du)\mu^{\delta}(dv)\}^2
P_{\lambda}\{T^{\delta}\neq\Delta \}.\eqno(2.14)$$
We finally arrive at
$$c_5\delta^2D^2(\varepsilon)
P_{\lambda}\left\{T^{\delta}\neq\Delta\right\}
+c_6\int_{\bbR^d}e^{-\frac{\varepsilon^2}{2}|\xi|^2}Q_{\mu^{\delta}}(\xi)d\xi$$
$$\ge c_7\left(\int_{\bbR^d}e^{-\frac{\varepsilon^2}{2}|\xi|^2}Q_{\mu^{\delta}}(\xi)d\xi\right)^2
P_{\lambda}\left\{T^{\delta}\neq\Delta\right\},\eqno(2.15)$$
where
$c_5,~c_6,~c_7\in(0,\infty)$ are some constants completely independent of
$\delta$ and $\varepsilon.$

There are four small steps from (2.14) to (2.15):
A. For any fixed $v\in\bbR^n_+,$
$\sum_{\pi}\int_{u\succeq_{\pi} v}(\bullet)\mu^{\delta}(du)=
\int_{\bbR^n_+}(\bullet)\mu^{\delta}(du).$
B. Use (2.12). C. There is no problem with interchanging the order of integration
once more thanks to the term $e^{-\frac{\varepsilon^2}{2}|\xi|^2}.$
D. Use (2.10).

\vskip .15in
Choose any sequence $\delta_k\downarrow0$ as $k\rightarrow\infty$ where $k=1,2,\cdots.$
Since $G^{\delta_1}$ is bounded, there exists a probability measure $\mu$
such that along some subsequence
$\delta_m\rightarrow0,$
$\mu^{\delta_m}\rightarrow\mu$ weakly.
To see that $\mu$ is supported on $G$, we notice that
$G$, as well as each $G^{\delta},$ is compact and that
$G\subset G^{\delta_{m+1}}\subset
G^{\delta_{m}}.$ Taking the indicator function $1_{G^{\delta}}$ and noting that 
$\mu^{\delta_m}$ is supported on $G^{\delta_{m}}$, we can easily find a contradiction
if $\mu$ has a positive mass on a compact set $B$ with $B\cap G=\emptyset.$ 
Next we write
$$\int_{\bbR^d}e^{-\frac{\varepsilon^2}{2}|\xi|^2}Q_{\mu^{\delta_m}}(\xi)d\xi
=\int_{\bbR^n_+}\int_{\bbR^n_+}f(s,t)\mu^{\delta_m}(ds)\mu^{\delta_m}(dt)$$
where
$$
f(s,t)=
\int_{\bbR^d}e^{-\frac{\varepsilon^2}{2}|\xi|^2}
e^{-\sum^n_{i=1}|t_i-s_i|\Psi'_i(\mbox{sgn}(t_i-s_i)\xi)}
\prod_{j=1}^N\mbox{Re}\left(\frac{1}{1+\Psi_j(\xi)}\right)d\xi.$$
Quite clearly, $f(s,t)$ is a bounded continuous function.
From the approximation argument from simple functions to
bounded continuous functions in the weak convergence for probability measures,
it also holds that 
$\mu^{\delta_m}\rightarrow\mu$ weakly in the double space sense:
$$\int_{\bbR^n_+}\int_{\bbR^n_+}f(s,t)\mu^{\delta_m}(ds)\mu^{\delta_m}(dt)
\rightarrow
\int_{\bbR^n_+}\int_{\bbR^n_+}f(s,t)\mu(ds)\mu(dt).$$
In other words,
$$\lim_{m\rightarrow\infty}
\int_{\bbR^d}e^{-\frac{\varepsilon^2}{2}|\xi|^2}Q_{\mu^{\delta_m}}(\xi)d\xi
=\int_{\bbR^d}e^{-\frac{\varepsilon^2}{2}|\xi|^2}Q_{\mu}(\xi)d\xi>0.\eqno(2.16)$$
Now rewrite (2.15) as
$$c_5\delta^2D^2(\varepsilon)
P_{\lambda}\left\{T^{\delta}\neq\Delta\right\}
\left(\int_{\bbR^d}e^{-\frac{\varepsilon^2}{2}|\xi|^2}Q_{\mu^{\delta}}(\xi)d\xi\right)^{-2}$$
$$+c_6\left(\int_{\bbR^d}e^{-\frac{\varepsilon^2}{2}|\xi|^2}Q_{\mu^{\delta}}(\xi)d\xi\right)^{-1}
\ge c_7
P_{\lambda}\left\{T^{\delta}\neq\Delta\right\}.\eqno(2.17)$$

Recall that
$$P_{\lambda}\left\{T^{\delta_m}\neq\Delta\right\}\rightarrow
\prod_{i=1}^{2^N}
E\left\{
\lambda_d\left(\overline{Z^i(G\times[0,l]^N)}\right)\right\}
\in(0,\infty)
\eqno(2.18)$$
downwards as $m\rightarrow\infty.$

It follows from (2.16), (2.17) and (2.18) that
$$c_6\left(\int_{\bbR^d}e^{-\frac{\varepsilon^2}{2}|\xi|^2}Q_{\mu}(\xi)d\xi\right)^{-1}
\ge
c_7
\prod_{i=1}^{2^N}
E\left\{
\lambda_d\left(\overline{Z^i(G\times[0,l]^N)}\right)\right\}
>0.\eqno(2.19)$$
Finally, let
$\varepsilon\rightarrow0$ in (2.19) to finish.

\vskip .25in

\noindent
{\bf Case 2} $\lambda_n(G)=0.$~~~This is the major case because we are more interested in
the measure $\mu$ on a nontrivial set $G$ with $\lambda_n(G)=0$
satisfying (1.1).

Fix a point $q=(q_1,\cdots,q_n)\in G.$ For $\eta\in(0,1],$ let
$G_{\eta}=G\bigcup
([q_1,q_1+\eta]\times\cdots\times[q_n,q_n+\eta]).$ Choose a copy $Y'$ of $Y$ independent of
$Z$. Consider the deterministic $(d,d)$ additive L\'evy process
$$\zeta^{\eta}_t=(at_1,at_2,\cdots,at_d),~~~(t_1,t_2,\cdots,t_d)\in\bbR^d_+,$$
where $a=\eta^{-2n(2^N-1)/d}.$ Define the additive L\'evy process
$$Z^{\eta}=Y'+\zeta^{\eta}.$$
Now we replace the $Z$ in Case 1 by an even larger additive L\'evy process
$$Z=(Z^1,Z^2,\cdots,Z^{2^N},Z^{\eta}).$$
Observe that
$$E\{\lambda_d(Y'(G_{\eta})+\zeta^{\eta}([0,1]^d))\}\ge
E\{\lambda_d(\zeta^{\eta}([0,1]^d))\}=a^d=\eta^{-2n(2^N-1)}.$$
By Proposition 10.3 and Theorem 2.1 of Khoshnevisan, Xiao and Zhong [1], for all $i\ge2,$ 
$$E\{\lambda_d(Z^i(G_{\eta}\times[0,l]^N))\}\ge
c_8
\eta^{2n}
E\{\lambda_d(Z^1(G_{\eta}\times[0,l]^N))\},$$
where $c_8\in(0,\infty)$ is some constant totally independent
of $\eta.$
As in Case 1, we first compute
$$P_{\lambda}
\left\{Z^i(G_{\eta}^{\delta}\times[0,l]^N)\bigcap
B(0,\delta)\neq\emptyset,~1\le i\le2^N,~ 
Z^{\eta}(G_{\eta}^{\delta}\times[0,1]^d)\bigcap
B(0,\delta)\neq\emptyset
\right\}.$$
We have
\begin{eqnarray*}
&&
P_{\lambda}
\left\{Z^i(G_{\eta}^{\delta}\times[0,l]^N)\bigcap
B(0,\delta)\neq\emptyset,~1\le i\le2^N,~ 
Z^{\eta}(G_{\eta}^{\delta}\times[0,1]^d)\bigcap
B(0,\delta)\neq\emptyset
\right\}\\
&&
=\prod_{i=1}^{2^N}
E\left\{
\lambda_d(Z^i(G_{\eta}^{\delta}\times[0,l]^N)+
B(0,\delta))\right\}\cdot
E\{\lambda_d(Y'(G_{\eta}^{\delta})+\zeta^{\eta}([0,1]^d)+B(0,\delta))\}\\
&&\ge
c_8^{2^N-1}
\left[E\{
\lambda_d(Z^1(G_{\eta}\times[0,l]^N))\}\right]^{2^N}
(\eta^{2n})^{2^N-1}
\eta^{-2n(2^N-1)}\\
&&
=c_8^{2^N-1}
\left[E\{
\lambda_d(Z^1(G_{\eta}\times[0,l]^N))\}\right]^{2^N}\\
&&\ge
c_8^{2^N-1}
\left[E\{\lambda_d(X([0,l]^N)+Y(G))\}\right]^{2^N}>0.
\end{eqnarray*}

Similarly, we define the random variable
$T^{\delta}$ as 
$T^{\delta}\neq\Delta$ if and only if 
$T^{\delta}\in{\bf Q}^n_+\bigcap(0,\infty)^n\bigcap G_{\eta}^{\delta}$
and
there exist
$t_i\in (0,l]^N\bigcap
{\bf Q}^N_+,$
$1\le i\le 2^N,$ such that
$|Z^{i}_{(T^{\delta},t_i)}|\le\delta$
and
there exists a
$t_0\in (0,1]^d\bigcap
{\bf Q}^d_+$
such that
$|Z^{\eta}_{(T^{\delta},t_0)}|\le\delta.$

We then redo the
$\mathcal{M}^{{\bf A},{\bf f},\nu}_{{\bf s}}$ thing as in Case 1, but
just do not do it for
$Z^{\eta}.$ We wind up with an inequality in (2.19):
$$\left(\int_{\bbR^d}e^{-\frac{\varepsilon^2}{2}|\xi|^2}Q_{\mu^{\eta}}(\xi)d\xi\right)^{-1}
\ge
c_9
\left[E\{\lambda_d(X([0,l]^N)+Y(G))\}\right]^{2^N}>0,\eqno(2.20)$$
where $c_9\in(0,\infty)$ is some constant completely independent of
$\eta$ and $\varepsilon,$
and $\mu^{\eta}$ is a probability measure on 
$G_{\eta}.$ By Fatou's lemma, we can find a probability measure $\mu$ on $G$ such
that
$$\left(\int_{\bbR^d}e^{-\frac{\varepsilon^2}{2}|\xi|^2}Q_{\mu}(\xi)d\xi\right)^{-1}
\ge
c_9
\left[E\{\lambda_d(X([0,l]^N)+Y(G))\}\right]^{2^N}.\eqno(2.21)$$
This time,
let
$\varepsilon\rightarrow0$ in (2.21) to complete the proof.
~~~~~~~~~~~~~~~~~~~~~~~~~~~~~~~~~~~~~~~~~~~~
~~~~~~~~~~~~$\Box$

\vskip 0.25in
\noindent
{\bf Proof of (1.3)}~~
In Theorem 1.1, we let $n=d$ and
$Y=\zeta$, where $\zeta_t=(t_1,t_2,\cdots,t_d),~(t_1,t_2,\cdots,t_d)\in\bbR^d_+.$
$\zeta_t$ is a deterministic additive L\'evy process with
L\'evy exponents $\Psi'_k(\xi)=-i\xi_k,~k=1,2,\cdots,d,$
$\xi=(\xi_1,\xi_2,\cdots,\xi_d)\in\bbR^d.$
First, assume that $F\subset\bbR^d_+.$
Note that $\zeta(F)=F.$ 
The reader can check that for any probability measure $\mu$ in $\bbR^d_+$, 
$$|\hat{\mu}(\xi)|^2=
\int_{\bbR^d_+}\int_{\bbR^d_+}
e^{-\sum^d_{k=1}|t_k-s_k|\Psi'_k(\mbox{sgn}(t_k-s_k)\xi)}\mu(ds)\mu(dt),~~~\xi\in\bbR^d.$$
(1.3) follows in the case when $F\subset\bbR^d_+$.

\vskip .15in

Next we consider an arbitrary $F$. Let $\bbR^d_p$ be the $p$-th closed quadrant of 
$\bbR^d.$ Define
$F^p=F\bigcap\bbR^d_p.$ Since 
$X(\bbR^N_+)+F=\bigcup_p(X(\bbR^N_+)+F^p),$
$\lambda_d(X(\bbR^N_+)+F)\le\sum_p
\lambda_d(X(\bbR^N_+)+F^p).$ 
Thus,
$$E\{\lambda_d(X(\bbR^N_+)+F)\}>0\Longleftrightarrow
E\{\lambda_d(X(\bbR^N_+)+F^p)\}>0$$
for some $p$.

Let $\zeta^p$ be the deterministic additive L\'evy process corresponding to
the quadrant $\bbR^d_p;$
i.e.,
$\zeta^p_t=(\pm t_1,\pm t_2,\cdots,\pm t_d),~(t_1,t_2,\cdots,t_d)\in\bbR^d_+.$
Let $(\Psi'_1,\cdots, \Psi'_d)$ be the L\'evy exponent of
$\zeta^p$. Let $\widetilde{F}^p=(\zeta^p)^{-1}(F^p),$ where
$(\zeta^p)^{-1}:\bbR^d_p\rightarrow\bbR^d_+$ is the inverse of $\zeta^p$.
Then  
$\widetilde{F}^p\subset\bbR^d_+$ and
$\zeta^p(\widetilde{F}^p)=F^p.$
Let $\tilde{\mu}$ be any probability measure in
$\bbR^d_+.$ Then $\mu=\tilde{\mu}\circ(\zeta^p)^{-1}$
is a probability measure in 
$\bbR^d_p$. In particular, if $\tilde{\mu}$ is on $\widetilde{F}^p$, then
$\mu$ is on $F^p.$
Similarly, we have
$$|\hat{\mu}(\xi)|^2=
\int_{\bbR^d_+}\int_{\bbR^d_+}
e^{-\sum^d_{k=1}|t_k-s_k|\Psi'_k(\mbox{sgn}(t_k-s_k)\xi)}
\tilde{\mu}(ds)\tilde{\mu}(dt),~~~\xi\in\bbR^d.\eqno\Box$$

\vskip .3in

\centerline{REFERENCES}
\vskip .2in

\noindent {[1]} Khoshnevisan, D., Xiao, Y. and Zhong, Y. (2003). 
Measuring the range of an additive 

\indent \indent
L\'evy process. {\it Ann. Probab.} {\bf 31}, 1097-1141.

\noindent {[2]} Khoshnevisan, D. and Xiao, Y. (2005). 
L\'evy processes: capacity and Hausdorff  
dimension.

\indent \indent
{\it Ann. Probab.} {\bf 33}, 841-878.

\end{document}